# A note on Cardano's formula

Dietmar Pfeifer

Institute of Mathematics, Carl von Ossietzky Universität Oldenburg, D-26111 Oldenburg, Germany

**Abstract**

In this paper we show that Cardano's formula for the solution of cubic equations can be reduced to expressions involving only square roots of rational numbes if the real root itself is rational.

## 1. Introduction and main result

By 1494, cubic equations were in general unsolvable algebraically. The first abstract solutions can be found in the work of the Italian mathematicians Scipione del Ferro, Niccolò Tartaglia, Girolamo Cardano, and Rafael Bombelli (cf. [1], Chapter 12: Algebra in the Renaissance, p.383ff). Since Cardano was obviously the first to publish the corresponding results in 1545, the formulas for the solution found in that time are named after him.

Another simple approach to an algebraic solution of cubic equations can be found in [3], p. 134 ff, analyzing the divisors of the integer constant term in the equation motivated by Vieta's theorem. A completely different way to solve cubic equations "algebraically" can be found in [2], p. 235, which, however, turns out to be only a very good numerical approximation of the solution.

To be more precise, we consider the equation $x^3 + 3ax = 2b$, $x \in \mathbb{C}$ for the case $a^3 + b^2 \geq 0$ with rational $a,b$. This is the case where exactly one root is real and the other ones are complex. According to Cardano's formula (cf. [2], p. 233) we can write the real solution $x$ as

$$x = w_1 - w_2 \text{ where } w_1 = \sqrt[3]{\sqrt{a^3 + b^2} + b}, \; w_2 = \sqrt[3]{\sqrt{a^3 + b^2} - b}. \tag{1}$$

However, even if the solution $x$ is rational, $w_1$ and $w_2$ are typically irrational. We want to show in this note that in this case $w_1$ and $w_2$ can surprisingly be expressed as terms solely depending on square roots of rational numbers, i.e.

$$x = w_3 - w_4 \text{ where } w_3 = s\sqrt{a^3 + b^2} + t = w_1, \; w_4 = s\sqrt{a^3 + b^2} - t = w_2 \tag{2}$$

with $t = \dfrac{x}{2}$, $s = \dfrac{t}{b - 2at}$. This problem has also been addressed in [4] p. 163ff. and [5], p.140ff. On the contrary, our paper presents a complete explicit solution for all relevant cases.

**Example 1.** Choose $a = \dfrac{4}{3}$ and $b = \dfrac{75}{16}$. Then the resulting equation is $x^3 + 4x = \dfrac{75}{8}$ with the only real solution $x = \dfrac{3}{2}$. According to Cardano's formula, we have

$$w_1 = \sqrt[3]{\frac{43}{144}\sqrt{273} + \frac{75}{16}} = 2{,}126892637 \text{ and } w_2 = \sqrt[3]{\frac{43}{144}\sqrt{273} - \frac{75}{16}} = 0{,}626892636 \qquad (3)$$

with the numerical approximation $x = w_1 - w_2 = 1{,}500000001$. With $t = \frac{3}{4}$ and $s = \frac{12}{43}$ we obtain

$$w_3 = \frac{1}{12}\sqrt{273} + \frac{3}{4} = 2{,}126892637 \text{ and } w_4 = \frac{1}{12}\sqrt{273} - \frac{3}{4} = 0{,}626892637 \qquad (4)$$

with the exact solution $x = w_3 - w_4 = \frac{3}{2}$.

**2. Proof** of alternate solution form.

We write for short $r = w_2$. Then $w_1 = x + r$ and

$$\sqrt{a^3 + b^2} + b = (x+r)^3 = x^3 + 3x^2 r + 3xr^2 + r^3$$
$$= x^3 + 3x^2 r + 3xr^2 + \sqrt{a^3 + b^2} - b, \qquad (5)$$

hence

$$x^3 + 3x^2 r + 3xr^2 = 2b = x^3 + 3ax \quad \text{or} \quad r^2 + xr = a \qquad (6)$$

with the solution (note $r \geq 0$)

$$r = \frac{1}{2}\sqrt{x^2 + 4a} - \frac{x}{2}. \qquad (7)$$

We put $t := \frac{x}{2}$ and equate the above equations to obtain

$$\left(s\sqrt{a^3+b^2} + t\right)^3 = \left(s^3(a^3+b^2) + 3st^2\right)\sqrt{a^3+b^2} + 3s^2 t(a^3+b^2) + t^3 = \sqrt{a^3+b^2} + b \qquad (8)$$

and

$$\left(s\sqrt{a^3+b^2} - t\right)^3 = \left(s^3(a^3+b^2) + 3st^2\right)\sqrt{a^3+b^2} - 3s^2 t(a^3+b^2) - t^3 = \sqrt{a^3+b^2} - b. \qquad (9)$$

This means that we have to solve

$$t^3 + 3s^2 t(a^3+b^2) = b, \quad s^3(a^3+b^2) + 3st^2 = 1. \qquad (10)$$

Rewriting this as

$$t^2 + 3s^2(a^3+b^2) = \frac{b}{t}, \quad 3s^2(a^3+b^2) + 9t^2 = \frac{3}{s} \qquad (11)$$

we get by subtraction

$$8t^2 = \frac{3}{s} - \frac{b}{t} \text{ or } 0 = 8t^3 - \frac{3t}{s} + b = x^3 - \frac{3x}{2s} + b = 3ax - \frac{3}{2s}x + 3b \text{ or } s = \frac{x}{2(b-ax)} = \frac{t}{b-2at}. \qquad (12)$$

**Example 2.** Consider the equation $x^3 + 3ax = a^2 - a = 2b$ with the only real solution $x = \sqrt[3]{a^2} - \sqrt[3]{a}$. This can in general not be expressed solely with square roots of rational numbers.

Note that the statements above are not only restricted to Cardano's case but apply also to the *casus irreducibilis* where $a^3 + b^2 < 0$ with rational $a,b$. In this case we have only real solutions in spite of the fact that the roots in Cardano's formula are complex. We only have to keep in mind that $w_1$ and $w_2$ have in general three different representations as complex numbers. If $w_1$ and $w_2$ are any given complex values then the other two are obtained by multiplications with the two non-trivial unit roots $\varepsilon_1 = -\frac{1}{2} + \frac{1}{2}\sqrt{3}i$ and $\varepsilon_2 = -\frac{1}{2} - \frac{1}{2}\sqrt{3}i$.

**Example 3.** Consider the equation $x^3 - x = -\frac{3}{8}$ with $a = -\frac{1}{3}$ and $b = -\frac{3}{16}$. Note that $x = \frac{1}{2}$ is a rational solution to the cubic equation. Cardano's formula gives

$$w_1 = \sqrt[3]{-\frac{3}{16} + \frac{\sqrt{39}}{144}i} \text{ and } w_2 = \sqrt[3]{\frac{3}{16} + \frac{\sqrt{39}}{144}i} \qquad (13)$$

with the three corresponding numerical complex representations

$$w_1 \in \begin{Bmatrix} -0,5756939096 + 0,04370189911\,i, \\ 0,2500000001 - 0,5204165000\,i, \\ 0,3256939095 + 0,4767146009\,i \end{Bmatrix} \text{ and } w_2 \in \begin{Bmatrix} 0,5756939095 + 0,4370189898\,i, \\ -0,2500000000 - 0,5204165000\,i, \\ -0,3256939095 + 0,04767146011\,i \end{Bmatrix} \quad (14)$$

with the three numerical approximations

$x \in \{-1,151387819,\ 0,5000000001,\ 0,6513878189\}$. Our alternative approach gives $t = \frac{1}{4}$ and $s = -12$ with

$$w_3 = \frac{1}{4} - \frac{\sqrt{39}}{12}i = 0,25 - 0,5204165000\,i \text{ and } w_4 = -\frac{1}{4} - \frac{\sqrt{39}}{12}i = 0,25 - 0,5204165000\,i$$

with the exact solution $x = \frac{1}{2}$.

**Final remark.** Note that although in case of a rational solution to the cubic equation, relation (2) above does not provide an algorithmic approach to find the root of the equation. Note also that relation (2) remains true in the general case, for non-rational solutions of the cubic equation.

**Example 4.** Chose $a = 2$, $b = a^2 - a = 2$ (cf. Example 3). Then by a numerical approach, we obtain

$x = 0,6258168190$ with

$$w_1 = \sqrt[3]{\sqrt{a^3 + b^2} + b} = 1,761325364,\ w_2 = \sqrt[3]{\sqrt{a^3 + b^2} - b} = 1,135508545 \qquad (15)$$

according to relation (1) and

$$t = \frac{x}{2} = 0{,}3129084095, \quad s = \frac{t}{b-2at} = 0{,}4181219592, \qquad (16)$$

$$w_3 = s\sqrt{a^3+b^2} + t = 1{,}761325364 = w_1, \quad w_4 = s\sqrt{a^3+b^2} - t = 1{,}135508545 = w_2. \qquad (17)$$